\documentclass[reqno, 11pt]{amsart}
\usepackage{txfonts}
\usepackage{bm}
\usepackage{amsmath}
\usepackage{amsfonts}
\usepackage{amssymb}
\usepackage{amsthm}
\usepackage{graphicx}
\usepackage{color}
\usepackage{enumitem}

\theoremstyle{plain}

\newtheorem{theorem}{Theorem}[section]
\newtheorem{corollary}[theorem]{Corollary}
\newtheorem{definition}[theorem]{Definition}
\newtheorem{example}[theorem]{Example}
\newtheorem{lemma}[theorem]{Lemma}
\newtheorem{proposition}[theorem]{Proposition}
\newtheorem{remark}[theorem]{Remark}

\numberwithin{equation}{section}

\DeclareMathAlphabet\scr{U}{scr}{m}{n}
\SetMathAlphabet\scr{bold}{U}{scr}{b}{n}
  \DeclareFontFamily{U}{scr}{\skewchar\font'177}%
  \DeclareFontShape{U}{scr}{m}{n}{<-6>rsfs5<6-8>rsfs7<8->rsfs10}{}%
  \DeclareFontShape{U}{scr}{b}{n}{<-6>rsfs5<6-8>rsfs7<8->rsfs10}{}%

\newcommand{\auf}{[\![}
\newcommand{\zu}{]\!]}
\newcommand{\mal}{\stackrel{\mbox{\tiny$\bullet$}}{}}

\newcommand{\rr}{\mathbb R}
\newcommand{\rp}{\mathbb R _+}

\newcommand{\F}{\scr F}

\newcommand{\E}{\scr E}

\newcommand{\B}{\scr B}

\newcommand{\til}{\widetilde}
\renewcommand{\Re}{\mathrm{Re}}

\renewcommand{\epsilon}{\varepsilon}
\renewcommand{\theta}{\vartheta}
\renewcommand{\rho}{\varrho}

\begin{document}
\title[]{A characterization of the martingale property of exponentially affine processes}
\author[E.~Mayerhofer]{Eberhard Mayerhofer}
\address{Vienna Institute of Finance, University of Vienna and Vienna University of Economics
and Business Administration, Heiligenst\"adterstrasse 46-48, 1190
Vienna, Austria} \email{eberhard.mayerhofer@vif.ac.at}
\author[J.~ Muhle-Karbe]{Johannes Muhle-Karbe}
\address{Departement Mathematik, ETH Z\"urich, R\"amistrasse 101, 8092 Z\"urich, Switzerland}
\email{johannes.muhle-karbe@math.ethz.ch}
\author[A.~G.~Smirnov]{Alexander G. Smirnov}
\address{I.~E.~Tamm Theory Department, P.~N.~Lebedev Physical Institute, Leninsky prospect 53, Moscow 119991, Russia}
\email{smirnov@lpi.ru}
\thanks{E.M. and A.G.S gratefully acknowledge financial support from WWTF (Vienna Science and Technology Fund).
The research of A.G.S. was also supported by the Russian Foundation
for Basic Research (Grant No.~09-01-00835) and the Program for
Supporting Leading Scientific Schools (Grant No.~LSS-1615.2008). J.
M.-K. gratefully acknowledges support from the FWF (Austrian Science
Fund) under grant P19456. The authors also thank two anonymous referees,
whose insightful reports led to a considerable improvement of the
present article.}

\begin{abstract}
We consider local martingales of exponential form $M=e^X$ or $\E(X)$
where $X$ denotes one component of a multivariate affine process in
the sense of Duffie, Filipovi\'c and Schachermayer
\cite{duffie.al.03}.  By completing the characterization of
conservative affine processes in \cite[Section 9]{duffie.al.03}, we
provide deterministic necessary and sufficient conditions in terms
of the parameters of $X$ for $M$ to be a true martingale.
\end{abstract}
 \maketitle
\section{Introduction}
A classical question in probability theory comprises the following.
Suppose the ordinary resp.\ stochastic exponential $M=\exp(X)$
resp.\ $\E(X)$\footnote{The \emph{stochastic exponential} $\E(X)$ of a semimartingale $X$ is the unique solution of the linear SDE $d\E(X)_t=\E(X)_{t-}dX_t$ with $\E(X)_0=1$, cf., e.g., \cite[I.4.61]{js.87} for more details.} of some process $X$ is a
positive \emph{local martingale} and hence a supermartingale. Then
under what (if any) additional assumptions is it in fact a
\emph{true martingale}?

This seemingly technical question is of considerable interest in
diverse applications, for example, absolute continuity of
distributions of stochastic processes (cf., e.g.,
\cite{cheridito.al.05} and the references therein), absence of
arbitrage in financial models (see, e.g.,
\cite{delbaen.schachermayer.95c})  or verification of optimality in
stochastic control (cf., e.g.,  \cite{elkaroui.81}).

In a general semimartingale setting it has been shown in
\cite{foellmer.72} that any supermartingale $M$ is a martingale if
and only if it is non-explosive under the associated \emph{F\"ollmer
measure} (also cf.\ \cite{yoerp.76}). However, this general result
is hard to apply in concrete models, since it is expressed in purely
probabilistic terms. Consequently, there has been extensive research
focused on exploiting the link between martingales and non-explosion
in various more specific settings, see, e.g., \cite{wong.heyde.04}.
In particular,  \emph{deterministic} necessary and sufficient
conditions for the martingale property of $M$ have been obtained if
$X$ is a one-dimensional diffusion (cf., e.g.,
\cite{delbaen.shirakawa.02, blei.engelbert.09} and the references
therein; also compare \cite{mijatovic.urusov.10}).

For processes with jumps, the literature is more limited and mostly
focused on sufficient criteria as in \cite{lepingle.memin.78,
kallsen.shiryaev.00b, protter.shimbo.08, kallsen.muhlekarbe.08b}. By
the independence of increments and the L\'evy-Khintchine formula, no
extra assumptions are needed for $M$ to be a true martingale if $X$
is a L\'evy process. For the more general class of \emph{affine
processes} characterized in \cite{duffie.al.03} the situation
becomes more involved. While no additional conditions are needed for
continuous affine processes, this no longer remains true in the
presence of jumps (cf.\ \cite[Example 3.11]{kallsen.muhlekarbe.08b}).
In this situation a necessary and sufficient condition for
one-factor models has been established in \cite[Theorem
2.5]{kellerressel.09}, whereas easy-to-check sufficient conditions
for the general case are provided by \cite[Theorem
3.1]{kallsen.muhlekarbe.08b}.

In the present study, we complement these results by sharpening
\cite[Theorem 3.1]{kallsen.muhlekarbe.08b} in order to provide
deterministic necessary and sufficient conditions for the martingale
property of $M=\E(X^i)$ resp.\ $\exp(X^i)$ in the case where $X^i$
is one component of a general non-explosive affine process $X$. As
in \cite{kellerressel.09,kallsen.muhlekarbe.08b} these conditions
are expressed in terms of the admissible \emph{parameters} which
characterize the distribution of $X$ (cf.\ \cite{duffie.al.03}).

Since we also use the linkage to non-explosion, we first complete
the characterization of \emph{conservative}, i.e.\ non-explosive,
affine processes from \cite[Section 9]{duffie.al.03}. Generalizing
the arguments from \cite{kallsen.muhlekarbe.08b}, we then establish
that $M$ is a true martingale if and only if it is a local
martingale and a related affine process is conservative. Combined
with the characterization of local martingales in terms of
semimartingale characteristics  \cite[Lemma 3.1]{kallsen.03} this
then yields necessary and sufficient conditions for the martingale
property of $M$.

The article is organized as follows. In Section \ref{se: prelim}, we recall terminology and results on affine Markov processes from \cite{duffie.al.03}. Afterwards, we characterize conservative affine processes. Subsequently, in Section \ref{sec: exp}, this characterization is used to provide necessary and sufficient conditions for the martingale property of exponentially affine processes. Appendix \ref{sec: ODEs} develops ODE comparison results in a general non-Lipschitz setting that are used to establish the results in Section \ref{sec: cons}.

\section{Affine processes}\label{se: prelim}
For stochastic background and terminology, we refer to
\cite{js.87,revuz.yor.99}. We work in the setup of
\cite{duffie.al.03}, that is we consider a time-homogeneous Markov
process with  state space $D:=\rp^m \times \rr^n$, where $m, n \geq
0$ and $d=m+n \geq 1$. We write $p_t(x,d\xi)$ for its transition
function and let $(X,\mathbb{P}_x)_{x \in D}$ denote its realization
on the canonical filtered space $(\Omega,\scr{F}^0,(\scr{F}^0_t)_{t
\in \rp})$ of paths $\omega: \rp  \to D_{\Delta}$ (the
one-point-compactification of $D$). For every $x \in D$,
$\mathbb{P}_x$ is a probability measure on $(\Omega,\F^0)$ such that
$\mathbb{P}_x(X_0=x)=1$ and the Markov property holds, i.e.\
\begin{eqnarray*}
\mathbb{E}_x(f(X_{t+s})|\F^0_s)&=&\int_{D} f(\xi)p_{t}(X_{s},d\xi)\\
&=&\mathbb{E}_{X_s}(f(X_t)), \quad \mathbb{P}_x \textrm{--a.s.}\
\quad \forall t,s,\in \rp,
\end{eqnarray*}
for all bounded Borel-measurable functions $f: D \to \mathbb{C}$.
The Markov process $(X,\mathbb{P}_x)_{x \in D}$ is called
\emph{conservative} if $p_t(x,D)=1$, \emph{stochastically
continuous} if we have $p_s(x,\cdot) \to p_t(x,\cdot)$ weakly on
$D$, for $s \to t$, for every $(t,x) \in \rp \times D$, and
\emph{affine} if, for every $(t,u) \in \rp \times i\rr^d$, the
characteristic function of $p_t(x,\cdot)$ is of the form
\begin{equation}\label{e:affine}
\int_D e^{\langle u,\xi
\rangle}p_t(x,d\xi)=\exp\left(\psi_0(t,u)+\langle \psi(t,u), x
\rangle\right), \quad \forall x \in D,
\end{equation}
for some $\psi_0(t,u) \in \mathbb{C}$ and
$\psi(t,u)=(\psi_1(t,u),\ldots,\psi_d(t,u)) \in \mathbb{C}^d$.
Note that $\psi(t,u)$ is uniquely specified by \eqref{e:affine}. But $\mathrm{Im}(\psi_0(t,u))$ is only determined up to multiples of $2\pi$. As usual in the literature, we enforce uniqueness by requiring the continuity of $u \mapsto \psi_0(t,u)$ as well as $\psi_0(t,0)=\log(p_t(0,D)) \in (-\infty,0]$ (cf., e.g., \cite[\S 26]{bauer.02}).

For every stochastically continuous affine process, the mappings
$(t,u) \mapsto \psi_0(t,u)$ and $(t,u) \mapsto \psi(t,u)$ can be
characterized in terms of the following quantities:

\begin{definition}\label{definition par}
Denote by $h=(h_1,\ldots,h_d)$ the truncation function on $\rr^d$
defined by
$$h_k(\xi):=\begin{cases} 0, &\mbox{if } \xi_k=0, \\ (1 \wedge |\xi_k|)\frac{\xi_k}{|\xi_k|}, &\mbox{otherwise.} \end{cases} $$
Parameters $(\alpha,\beta,\gamma,\kappa)$ are called
\emph{admissible}, if
\begin{itemize}
\item $\alpha=(\alpha_0,\alpha_1,\ldots,\alpha_d)$ with symmetric positive semi-definite $d \times d$-matrices $\alpha_j$ such that $\alpha_j=0$ for $j \geq m+1$ and $\alpha_j^{kl}=0$ for $0 \leq j \leq m$, $1 \leq k,l \leq m$ unless $k=l=j$;
\item $\kappa=(\kappa_0,\kappa_1,\ldots,\kappa_d)$ where $\kappa_j$ is a Borel measure on $D \backslash \{0\}$ such that $\kappa_j=0$ for $j \geq m+1$ as well as $\int_{D \backslash \{0\}} ||h(\xi)||^2 \kappa_j(d\xi)<\infty$ for $0 \leq j \leq m$ and
$$\int_{D \backslash \{0\}} \vert h_k(\xi) \vert\kappa_j(d\xi)<\infty, \quad 0 \leq j \leq m, \quad 1 \leq k \leq m, \quad k \neq j;$$
\item $\beta=(\beta_0,\beta_1,\ldots,\beta_d)$ with $\beta_j \in \rr^d$ such that $\beta_j^k=0$ for $j \geq m+1$, $1 \leq k \leq m$ and
$$ \beta_j^k -\int_{D \backslash \{0\}} h_k(\xi)\kappa_j(d\xi) \geq 0, \quad 0 \leq j \leq m, \quad 1 \leq k \leq m, \quad k \neq j.$$
\item $\gamma=(\gamma_0,\gamma_1,\ldots,\gamma_d)$, where $\gamma_j \in \rp$ and $\gamma_j=0$ for $j=m+1,\dots,d$.
\end{itemize}
\end{definition}

Affine Markov processes and admissible parameters are related as
follows (cf.\ \cite[Theorem 2.7]{duffie.al.03} and \cite[Theorem
5.1]{kellerressel.al.09}):

\begin{theorem}\label{t:2.7}
Let $(X,\mathbb{P}_x)_{x \in D}$ be a stochastically continuous
affine process. Then there exist admissible parameters
$(\alpha,\beta,\gamma,\kappa)$ such that $\psi_0(t,u)$ and
$\psi(t,u)$ are given as solutions to the \emph{generalized Riccati
equations}
\begin{align}\label{e:riccati}
\partial_t \psi(t,u)&=R(\psi(t,u)), \quad \,\,\,\psi(0,u)=u,\\\
\partial_t\psi_0(t,u)&=R_0(\psi(t,u)),\quad\psi_0(0,u)=0,\label{e:riccati2}
\end{align}
where $R=(R_1,\dots, R_d)$ and for $0 \leq i \leq d$,
\begin{equation}\label{e:R}
R_i(u):=\frac{1}{2}\langle \alpha_i u, u\rangle+\langle \beta_i, u
\rangle-\gamma_i+\int_{D \backslash \{0\}} \left(e^{\langle
u,\xi\rangle} -1-\langle u, h(\xi) \rangle\right)\kappa_i(d\xi).
\end{equation}
Conversely, for any set $(\alpha,\beta,\gamma,\kappa)$ of admissible
parameters there exists a unique stochastically continuous affine
process such that \eqref{e:affine} holds for all $(t,u) \in \rp
\times i\rr^d$, where $\psi_0$ and $\psi$ are given by
\eqref{e:riccati2} and \eqref{e:riccati}.
\end{theorem}

Since any stochastically continuous affine process
$(X,\mathbb{P}_x)_{x \in D}$ is a Feller process (cf.\ \cite[Theorem
2.7]{duffie.al.03}), it admits a c\`adl\`ag modification and hence
can be realized on the space of c\`adl\`ag paths $\omega: \rp \to
D_{\Delta}$. If $(X,\mathbb{P}_x)_{x \in D}$ is also conservative it
turns out to be a semimartingale in the usual sense and hence can be
realized on the \emph{Skorokhod space}
$(\mathbb{D}^d,\scr{D}^d,(\scr{D}^d_t)_{t \in \rp})$ of $D$- rather
than $D_{\Delta}$-valued c\`adl\`ag paths. Here, $\scr{D}^d_t=\bigcap_{s>t}\scr{D}^{0,d}_s$ for the filtration $(\scr{D}^{0,d}_t)_{t \in \mathbb{R}_+}$ generated by $X$. The semimartingale characteristics of $(X,\mathbb{P}_x)_{x \in D}$ are then given in terms of the admissible parameters:

\begin{theorem}\label{t:2.12}
Let $(X,\mathbb{P}_x)_{x \in D}$ be a conservative, stochastically
continuous affine process and let $(\alpha,\beta,\gamma,\kappa)$ be the
related admissible parameters. Then $\gamma=0$ and for any $x \in
D$, $X=(X^1,\ldots,X^d)$ is a semimartingale on
$(\mathbb{D}^d,\scr{D}^d,(\scr{D}^d_t)_{t \in \rp},\mathbb{P}_x)$
with characteristics $(B,C,\nu)$ given by
\begin{eqnarray}
B_t &=& \int_0^t \left(\beta_0+\sum_{j=1}^d \beta_j X^j_{s-}\right) ds,\label{e:b}\\
C_t &=& \int_0^t \left(\alpha_0+\sum_{j=1}^d \alpha_j X^j_{s-}\right) ds,\label{e:c}\\
\nu(dt,d\xi) &=& \left(\kappa_0(d\xi)+\sum_{j=1}^d X^j_{t-}
\kappa_j(d\xi)\right) dt,\label{e:nu}
\end{eqnarray}
relative to the truncation function $h$. Conversely, let $X'$ be a
$D$-valued  semimartingale defined on some filtered probability
space $(\Omega',\scr{F}',(\F'_t),\mathbb{P}')$. If
$\mathbb{P}'(X'_0=x)=1$ and $X'$ admits characteristics of the form
\eqref{e:b}-\eqref{e:nu} with $X_{-}$ replaced by $X_{-}'$, then
$\mathbb{P}' \circ X'^{-1}=\mathbb{P}_x$.
\end{theorem}

\begin{proof} $\gamma=0$ is shown in \cite[Proposition 9.1]{duffie.al.03}; the remaining assertions follow from \cite[Theorem 2.12]{duffie.al.03}.
\end{proof}

\section{Conservative affine processes}\label{sec: cons}
In view of Theorem \ref{t:2.12}, the powerful toolbox of
semimartingale calculus is made available for affine processes, provided that the Markov
process $(X,\mathbb P_x)_{x\in D}$ is
conservative. Hence, it is desirable to characterize this property
in terms of the parameters of $X$.  This is done in the present
section. The main result is Theorem \ref{th: char can cons},
which completes the discussion of conservativeness in \cite[Section
9]{duffie.al.03}.

To prove this statement, we proceed as follows. First, we recall some properties of the generalized Riccati equations \eqref{e:riccati}, \eqref{e:riccati2} established by Duffie et al.\ \cite{duffie.al.03}. In the crucial next step, we use the comparison results developed in the appendix to show that whereas the characteristic exponent $\psi$ of the affine process $X$ is not the \emph{unique} solution to these equations in general, it is necessarily the \emph{minimal} one among all such solutions. Using this observation, we can then show that conservativeness of the process $X$ is indeed \emph{equivalent} to uniqueness for the specific initial value zero. Note that \emph{sufficiency} of this uniqueness property was already observed in \cite[Proposition 9.1]{duffie.al.03}; here we show that this condition is also \emph{necessary}. 

Let us first introduce some definitions and notation. The partial
order on $\mathbb R^m$ induced by the natural cone $\mathbb R_+^m$
is denoted by $\preceqq$. That is, $x\preceqq 0$ if and only if
$x_i\leq 0$ for $i=1,\dots,m$. A function $g: D_g\rightarrow\mathbb
R^m$ is \emph{quasimonotone increasing} on $D_g\subset \mathbb R^m$
(\emph{qmi} in short, for a general definition see section \ref{sec:
ODEs}) if and only if for all $x,y\in D_g$ and $i=1,\dots,m$ the
following implication holds true:
\[
(x\preceqq y,\quad x_i=y_i)\;\Rightarrow \; g_i(x)\leq g_i(y).
\]
In the sequel we write $\mathbb R_{--}:=(-\infty,0)$ and $\mathbb
C_{--}:=\{c\in \mathbb C\,\mid\,\Re(c)\in \mathbb R_{--}\}$.
Moreover, we introduce the index set $\mathcal I:=\{1,\dots, m\}$
and, accordingly, define by $u_{\mathcal I}=(u_1,\dots, u_m)$ the
projection of the $d$--dimensional vector $u$ onto the first $m$
coordinates. Similarly $R_{\mathcal I}$ denotes
the first $m$ components of $R$, i.e. $R_\mathcal I=(R_1,\dots,R_m)$
and $R_\mathcal I(u_I,0):=(R_1(u_1,\dots,u_m,0,\dots,0), \dots,
R_m(u_1,\dots,u_m,0,\dots,0))$. Finally, $\psi_{\mathcal{I}}$ and $\psi_{\mathcal{I}}(t,(u_{\mathcal{I}},0))$ are defined analogously.

For this section the uniqueness of solutions to eqs.
\eqref{e:riccati}--\eqref{e:riccati2} is essential. It is adressed
in the following remark. For more detailed information, we refer to
\cite[Sections 5 and 6]{duffie.al.03}.
\begin{remark}\rm
\begin{enumerate}\label{rem: uniquenss}
\item \label{uniquenss issue 1} Due to the admissibility conditions
on the jump parameters $\kappa$ the domains of $R_0$ and $R$ can be be extended from $i\mathbb
R^d$ to $\mathbb C_{-}^m\times i\mathbb R^n$. Moreover, $R_0, R$ are analytic functions on $\mathbb C_{--}^m\times
i\mathbb R^n$, and admit a unique continuous extension to $\mathbb C_{-}^m\times
i\mathbb R^n$.
\item In general, $R$ is not locally Lipschitz on $i\mathbb R^d$,
but only continuous (see \cite[Example 9.3]{duffie.al.03}). This
lack of regularity prohibits to provide well-defined
$\psi_0,\psi$ by simply solving
\eqref{e:riccati}--\eqref{e:riccati2}, because unique
solutions do not always exist, again cf.\ \cite[Example 9.3]{duffie.al.03}. Hence another
approach to construct unique characteristic exponents $\psi_0,\psi$
is required. Duffie et al.\ \cite{duffie.al.03} tackle this problem by
first proving the existence of unique global solutions
$\psi_0^\circ,\psi^\circ$ on $\mathbb C_{--}^m\times i\mathbb R^n$,
where uniqueness is guaranteed by the analyticity of $R$, see
\ref{uniquenss issue 1}. Their unique continuous extensions to the
closure $\mathbb C_{-}^m\times i\mathbb R^n$ are also differentiable
and solve \eqref{e:riccati}--\eqref{e:riccati2} for $u\in i\mathbb
R^d$. Moreover, they satisfy  \eqref{e:affine}. Henceforth, $\psi_0,\psi$ therefore denote these unique extensions.
\end{enumerate}
\end{remark}

\begin{lemma}\label{lem: properties affine on canonical state}
The affine transform formula \eqref{e:affine} also holds for
$u=(u_{\mathcal I},0)\in \mathbb R_-^d$ with characteristic
exponents $\psi_0(t,(u_{\mathcal I},0)): \,\mathbb R_+\times \mathbb
R_-^m\rightarrow\mathbb R_-$ and $\psi_{\mathcal I}(t,(u_{\mathcal
I},0)): \,\mathbb R_+\times \mathbb R_-^m\rightarrow\mathbb R_-^m$
satisfying
\begin{align}\label{eq phi}
\partial_t{\psi}_0(t,(u_{\mathcal I},0))&=R_0((\psi_{\mathcal I}(t,(u_{\mathcal I},0)),0)),\quad\quad\,\,\,
\psi_0(0,(u_{\mathcal I},0))=0,\\\label{eq psi}\partial_t
{\psi_{\mathcal I}}(t,(u_\mathcal I,0))&=R_{\mathcal
I}((\psi_{\mathcal I}(t,(u_{\mathcal I},0)),0))\quad\quad\,\,\,
\psi_{\mathcal I}(0,(u_{\mathcal I},0))=u_{\mathcal I}.
\end{align}
Furthermore we have:
\begin{itemize}
\item $R_0, R_{\mathcal I}$ are continuous functions on $\mathbb
R_{-}^m$ such that $R_0(0)\leq 0$, $R_{\mathcal I}(0)\preceqq 0$
\item $R_{\mathcal I}((u_{\mathcal I},0))$ is locally Lipschitz continuous on
$\mathbb R_{--}^m$ and qmi on $\mathbb R_-^m$, \item $\psi_{\mathcal
I}(t,(u_{\mathcal I},0))$ restricts to
an $\mathbb{R}^m_{--}$-valued unique
global solution $\psi_{\mathcal I}^{\circ}(t,(u_{\mathcal I},0))$
of \eqref{eq psi} on $\mathbb R_+\times \mathbb R_{--}^m$.
\end{itemize}
\end{lemma}

\begin{proof}
By \cite{kellerressel.al.09}, any {\it stochastically continuous}
affine processes is {\it regular} in the sense of
\cite{duffie.al.03}. Hence, the first statement is a consequence of
\cite[Proposition 6.4]{duffie.al.03}. The regularity of $R_0$ and
$R_{\mathcal I}$ follows from \cite[Lemma 5.3 (i) and
(ii)]{duffie.al.03}. Equation \eqref{e:R} shows $R_0(0)\leq 0$ and
$R_{\mathcal I}(0)\preceqq 0$.  The mapping $v\mapsto R_{\mathcal
I}((v,0))$ is qmi on $\mathbb R_-^m$ by \cite[Lemma 4.6]{PhDMKR},
whereas the last assertion is stated in \cite[Proposition
6.1]{duffie.al.03}.
\end{proof}

In the following crucial step we establish the minimality of
$\psi_{\mathcal I}(t,(u_{\mathcal I},0))$ among all solutions of
\eqref{eq psi} with respect to the partial order $\preceqq$.
\begin{proposition}\label{prop: extremality}
Let $T>0$ and $u_{\mathcal{I}} \in \mathbb{R}_{-}^m$. If $g(t):
[0,T)\rightarrow \mathbb R_{-}^m$ is a solution of
\begin{equation}\label{e:eberhard1}
\partial_t g(t)=R_{\mathcal I}(g(t),0),
\end{equation}
subject to $g(0)=u_{\mathcal I}$, then $g(t)\succeqq \psi_{\mathcal
I}(t,(u_{\mathcal I},0))$, for all $t<T$.
\end{proposition}
\begin{proof}
The properties of $R_{\mathcal I}$ established in Lemma \ref{lem:
properties affine on canonical state} allow this conclusion by a use
of Corollary \ref{th1}. For an application of the latter, we make
the obvious choices $f=R_{\mathcal I}$, $D_f=\mathbb R_-^m$. Then we
know that for $u^\circ_{\mathcal I}\in \mathbb R_{--}^m$ we have
$g(t)\succeq \psi_{\mathcal I}^\circ(t,(u_{\mathcal I},0))$, for all
$t<T$. Now letting $u^\circ_{\mathcal I}\rightarrow u_{\mathcal I}$
and using the continuity of $\psi_{\mathcal I}$ as asserted in Lemma
\ref{lem: properties affine on canonical state} yields the
assertion.
\end{proof}

We now state the main result of this section, which is a full
characterization of conservative affine processes in terms of a
uniqueness criterium imposed on solutions of the corresponding
generalized Riccati equations. It is
motivated by a partial result of this kind provided in
\cite[Proposition 9.1]{duffie.al.03}, which gives a necessary
condition for conservativeness, as well as a sufficient one. Here, we show that their sufficient condition, which (modulo the assumption
$R(0)=0$) equals \ref{char1 point 3} below, is in fact also necessary for
conservativeness.
The proof is based on the comparison results for multivariate initial value problems developed in Appendix
\ref{sec: ODEs}.

\begin{theorem}\label{th: char can cons}
The following statements are equivalent:
\begin{enumerate}
\item \label{char1 point 1} $(X, \mathbb P_x)_{x\in D}$ is conservative,
\item \label{char1 point 3} $R_0(0)=0$ and there exists no non-trivial $ \mathbb R_-^m$-valued local solution $g(t)$ of \eqref{e:eberhard1} with $g(0)=0$.
\end{enumerate}
Moreover, each of these statements implies that $R(0)=0$.
\end{theorem}
\begin{proof}
\ref{char1 point 1}$\Rightarrow$\ref{char1 point 3}: By definition,
$X$ is conservative if and only if, for all $t\geq 0$ and $x\in D$, we
have
\[
1=p_t(x, D)=e^{\psi_0(t,0)+\langle
\psi(t,0),x\rangle}=e^{\psi_0(t,0)+\langle \psi_{\mathcal
I}(t,0),x_{\mathcal I}\rangle},
\]
because $\psi_i(t,(u_{\mathcal I},0))=0$, for $i=m+1,\dots,d$. By first putting $x=0$ and then using the arbitrariness of $x$, it follows that this is equivalent to
\begin{equation}\label{eq: ptD}
\psi_0(t,0)=0\textrm{  and  }\psi_{\mathcal
I}(t,0)=0, \quad \forall t \geq 0.
\end{equation}
Let $g$ be a (local) solution of \eqref{e:eberhard1} on some
interval $[0,T)$, satisfying $g(0)=0$ and with values in $\mathbb
R_{-}^m$. By Proposition \ref{prop: extremality}, $\psi(t,0)\preceqq
g(t)$, $0\leq t< T$. In view of \ref{char1 point 1} and
eq.~\eqref{eq: ptD}, the left side of the inequality is equal to
zero. This yields $g\equiv 0$. Now by Lemma \ref{lem: properties
affine on canonical state} and \ref{char1 point 1} (see \eqref{eq:
ptD})
\[
0=\psi_\mathcal I(t,0)=\int_0^t R_0(\psi_\mathcal I(s,0))ds=\gamma_0
t,\quad t\in [0,T),
\]
which implies $\gamma_0=0$ and hence \ref{char1 point
3}.

\ref{char1 point 3} $\Rightarrow$\ref{char1 point 1}: By Lemma
\ref{lem: properties affine on canonical state}, $g:=\psi_{\mathcal
I}(\cdot,0)$ is a solution of \eqref{e:eberhard1} with $g(0)=0$ and
values in $\mathbb{R}^m_{-}$. Assumption \ref{char1 point 3} implies
$\psi_{\mathcal I}(\cdot,0)\equiv 0$. Now $\gamma_0=R_0(0)=0$ as
well as $\psi_0(t_0,0)=0$ and \eqref{eq phi} yield $\psi_0(\cdot,
0)\equiv 0$. Hence \eqref{eq: ptD} holds and \ref{char1 point 1}
follows.

Finally, we show that either \ref{char1 point 1} or \ref{char1 point
3} implies $(\gamma_1,\dots,\gamma_m)=0$. Note that by
Definition \ref{definition par} we have
$\gamma_{m+1}=\dots=\gamma_d=0$. From \eqref{eq psi} for $u_\mathcal
I=0$ and from \eqref{eq: ptD} it follows that $0=R_j(0)\cdot t$ and
hence $R_j(0)=\gamma_j=0$ for all $1\leq j\leq m$.
\end{proof}

\begin{remark}\label{consremark}\rm
\begin{enumerate}
\item By Definition \ref{definition par}, $R_0(0)=0$, $R(0)=0$ is equivalent to $\gamma=0$. This means that the infinitesimal generator of the associated Markovian semi-group has zero potential, see \cite[Equation (2.12)]{duffie.al.03}. If an affine process with $\gamma=0$ fails to be conservative, then it must have state-dependent jumps.

\item \label{consremark for ODE result} The comparison results established in Appendix \ref{sec: ODEs} are the major tool for proving Proposition
  \ref{prop: extremality}. They are quite general and therefore allow
for a similar characterization of conservativeness of affine
processes on geometrically more involved state-spaces (as long as
they are proper closed convex cones). In particular, such a
characterization can be derived for affine processes on the cone of
symmetric positive semidefinite matrices of arbitrary dimension, see
also \cite[Remark 2.5]{CFMT}.
\item \label{expremark2cons} Conservativeness of $(X,\mathbb{P}_x)_{x \in D}$ and uniqueness for solutions of the ODE \eqref{e:eberhard1} can be ensured by requiring
\begin{equation}\label{e:sufficient}
\int_{D \backslash \{0\}} \left(|\xi_k| \wedge |\xi_k|^2\right)
\kappa_j(d\xi)<\infty, \quad 1 \leq k, j \leq m,
\end{equation}
as in \cite[Lemma 9.2]{duffie.al.03}, which implies that
$R_{\mathcal I}(\cdot,0)$ is locally Lipschitz continuous on
$\rr^m_{-}$.
\item \label{expremark3cons} If $m=1$, conservativeness corresponds to uniqueness of a one dimensional ODE and can be characterized more explicitly: \cite[Corollary 2.9]{duffie.al.03}, \cite[Theorem 4.11]{filipovic.01} and Theorem \ref{th: char can cons} yield that $(X,\mathbb{P}_x)_{x \in D}$ is conservative if and only if either \eqref{e:sufficient} holds or
    \begin{equation}\label{eq osgood}
\int_{0-} \frac{1}{R_{1}(u_1,0)}du_1=-\infty,
\end{equation}
where $\int_{0-}$ denotes an integral over an arbitrarily small left
neighborhood of $0$.
\end{enumerate}
\end{remark}

The sufficient condition \eqref{e:sufficient} from \cite[Lemma 9.2]{duffie.al.03} is easy to check in applications, since it can be read off directly from the parameters of $X$. However, the following example shows that it is not necessarily satisfied for conservative affine processes. 
This example is somewhat artificial and constructed so that the moment condition \eqref{e:sufficient} fails but the well-known Osgood condition \eqref{eq osgood} does not. While it is possible to extend the example in several directions (infinite activity, stable-like tails instead of discrete support, multivariate processes, etc.), we chose to present the simplest version in order to highlight the idea.

\begin{example}\label{crucial cons example}\rm
Define the measure
\[
\mu:= \sum_{n=1}^\infty \frac{\delta_n}{n^2},
\]
where $\delta_n$ is the Dirac measure supported by the one-point set
$\{n\}$. Then we have
\[
\beta_1:=\int_0^\infty h(\xi)\,d\mu(\xi)=\sum_{n=1}^\infty
\frac{1}{n^2}<\infty.
\]
Therefore the parameters $(\alpha,\beta,\gamma,\kappa)$ defined by
\[
\alpha=(0,0),\quad\beta=(0,\beta_1),
\quad\gamma=(0,0),\quad\kappa=(0,\mu)
\]
are admissible in the sense of Definition \ref{definition
par}.
Denote by $(X,\mathbb
P_x)_{x\in \mathbb R_+}$ the corresponding affine process provided
by Theorem \ref{t:2.7}. Then
\[
\int_0^\infty (\vert \xi \vert\wedge\vert
\xi\vert^2)\,d\mu(\xi)=\int_1^\infty \xi\,d\mu(\xi) =
\sum_{n=1}^\infty \frac{1}{n} = \infty,
\]
which violates the sufficient condition \eqref{e:sufficient} for
conservativeness. However, we now show that the necessary and sufficient condition \ref{char1 point 3}
of Theorem \ref{th: char can cons} is fulfilled, which in
turn ensures the conservativeness of $(X,\mathbb P_x)_{x\in \mathbb
R_+}$. By construction, $R_0(u)=0$ and \begin{equation}\label{1}
R(u)=R_1(u) = \int_1^\infty (e^{u\xi}-1)\,d\mu(\xi) =
\sum_{n=1}^\infty \frac{e^{un}-1}{n^2}.
\end{equation}
Clearly, $R(u)$ is smooth on $(-\infty,0)$, and differentiation of
the series on the right-hand side of~(\ref{1}) yields
\begin{align}
&R'(u)=\sum_{n=1}^\infty\frac{e^{un}}{n} \label{2} \\
&R''(u)=\sum_{n=1}^\infty e^{un} = \frac{e^{u}}{1-e^{u}}.\label{3}
\end{align}
By \eqref{3}, we have $R'(u)= -\ln(1-e^{u})+C$ and further by
\eqref{2}, $R'(u)$ tends to zero as $u\to-\infty$ and, therefore,
$C=0$. We thus obtain
\begin{equation}\label{4}
R'(u)= -\ln(1-e^{u}).
\end{equation}
Since $1-e^{u}=-u+O(u^2)$, we have $1-e^{u}\geq -u/2$ for $u\leq 0$
small enough. Hence,
\[
0 \leq R'(u)\leq -\ln \left(-\frac{u}{2}\right)
\]
for $u\leq 0$ small enough. As $R(0)=0$ by \eqref{1}, it follows
that
\begin{equation}\label{R est cons}
0 \geq R(u)=-\int_u^0 R'(u')\,du'\geq \int_u^0 \ln
\left(\frac{-u'}{2}\right)\,du' = -u\ln
\left(\frac{-u}{2}\right)+u\geq -2u\ln \left(\frac{-u}{2}\right)
\end{equation}
for $u\leq 0$ small enough. This implies
\begin{equation*}\label{int inf} \int_{-1}^{0-}\frac{du}{R(u)}=-\infty;
\end{equation*}
hence $(X,\mathbb P_x)_{x\in \mathbb R_+}$ is conservative by Remark
\ref{consremark} \ref{expremark3cons}.
\end{example}

\section{Exponentially affine martingales}\label{sec: exp}
We now turn to the characterization of exponentially affine
martingales. Henceforth, let $(X,\mathbb{P}_x)_{x\in D}$ be the
canonical realization on $(\mathbb{D}^d,\scr{D}^d,(\scr{D}^d_t)_{t
\in \rp})$ of a conservative, stochastically continuous affine
process with corresponding admissible parameters
$(\alpha,\beta,0,\kappa)$. 

We proceed as follows. First, we characterize the \emph{local} martingale property and the positivity of stochastic exponentials. Since these are ``local'' properties, they can be read directly from the parameters of the process. Afterwards, we consider the \emph{true} martingale property of $\E(X^i)$. Using Girsanov's theorem, we first establish that it is \emph{necessary} that a related affine process is conservative. Afterwards, we adapt the arguments from \cite{kallsen.muhlekarbe.08b} to show that this is also a \emph{sufficient} condition. Combined with the results of Section 3, this then characterizes the true martingale property of $\E(X^i)$ in terms of uniqueness of the solution of a system of generalized Riccati equations. Finally, we adapt our Example \ref{crucial cons example} to construct an exponentially affine local martingale $\E(X^i)$ for which the sufficient condition of \cite{kallsen.muhlekarbe.08b} fails, but uniqueness of the Riccati equations and hence the true martingale property of $\E(X^i)$ is assured by the Osgood condition \eqref{eq osgood}.

We begin with the local properties. Our first lemma shows that it can be read directly from the corresponding parameters whether $\E(X^i)$ is a local martingale.

\begin{lemma}\label{l:mloc}
Let $i \in \{1,\ldots,d\}$. Then $\E(X^i)$ is a local
$\mathbb{P}_x$-martingale for all $x \in D$ if and only if
\begin{equation}\label{e:integrable}
\int_{\{|\xi_i|>1\}}|\xi_i|\kappa_j(d\xi)<\infty, \quad 0 \leq j
\leq d,
\end{equation}
and
\begin{equation}\label{e:drift}
\beta_j^i +\int_{D \backslash \{0\}}
(\xi_i-h_i(\xi))\kappa_j(d\xi)=0, \quad 0 \leq j \leq d.
\end{equation}
\end{lemma}

\begin{proof} $\Leftarrow$: On any finite interval $[0,T]$, the mapping $t \mapsto X_{t-}$ is $\mathbb{P}_x$-a.s.\ bounded for all $x \in D$. Hence it follows from  Theorem \ref{t:2.12} and \cite[Lemma 3.1]{kallsen.03} that $X^i$ is a local $\mathbb{P}_x$-martingale. Since $\E(X^i)=1+\E(X^i)_{-} \mal X^i$ by definition of the stochastic exponential, the assertion now follows from \cite[I.4.34]{js.87}, because $\E(X^i)_{-}$ is locally bounded.\\
$\Rightarrow$: As $\kappa_j=0$ for $j=m+1,\ldots, d$ and $X^j_{-}$
is nonnegative for $j=1,\ldots,m$, \cite[Lemma 3.1]{kallsen.03} and
Theorem \ref{t:2.12} yield that
$\int_{\{|\xi_i|>1\}}|\xi_i|\kappa_0(d\xi)<\infty$ and
\begin{equation}\label{e:component}
\int_{\{|\xi_i|>1\}}|\xi_i|\kappa_j(d\xi)X^j_{-}<\infty, \quad 1
\leq j \leq m,
\end{equation}
up to a $d\mathbb{P}_x \otimes dt$-null set on $\Omega \times \rp$
for any $x \in D$. Now observe that \eqref{e:component} remains
valid if $X_{-}$ is replaced by $X$, because $X_{-}=X$ holds
$d\mathbb{P}_x \otimes dt$-a.e., for any $x \in D$. Setting
$\Omega_x=\{X_0=x\}$ for some $x \in D$ with $x^j>0$, the
right-continuity of $X$ shows that there exist $\epsilon>0$ and a
strictly positive random variable $\tau$ such that $X^j_{t}(\omega)
\geq \epsilon$ for all $0 \leq t \leq \tau(\omega)$ and for all
$\omega \in \Omega_x$. Denoting the set on which \eqref{e:component}
holds by $\til{\Omega}_0$,  it follows that the set $\til{\Omega}_0
\cap \auf 0,\tau \zu \cap  \Omega_x \times\mathbb{R}_+  \subset
\Omega \times \mathbb{R}_+$ has strictly positive $d\mathbb{P}_x
\otimes dt$-measure. Therefore it contains at least one $(\omega,t)$
for which
$$ \epsilon \int_{\{|\xi_i|>1\}} |\xi_i| \kappa_j(d\xi) \leq  \int_{\{|\xi_i|>1\}} |\xi_i|\kappa_j(d\xi) X^j_t(\omega) < \infty.$$
Hence \eqref{e:integrable} holds. We now turn to \eqref{e:drift},
which is well-defined by \eqref{e:integrable}. Set
$$\til{\beta}^i_j:= \beta^i_j +\int_{D \backslash \{0\}} (\xi_i-h_i(\xi)) \kappa_j(d\xi), \quad 0 \leq j \leq d.$$
Again by \cite[Lemma 3.1]{kallsen.03} and Theorem \ref{t:2.12}, we
have
\begin{equation}\label{e:componentdrift}
\til{\beta}^i_0+\sum_{j=1}^d \til{\beta}^i_j X^{j}_{-}=0,
\end{equation}
up to a $d\mathbb{P}_x \otimes dt$-null set on $\Omega \times \rp$
for all $x \in D$. As above, \eqref{e:componentdrift} remains valid
if $X_{-}$ is replaced by $X$. But now, using Fubini's theorem and
the right-continuity of $X$ we find that \eqref{e:componentdrift}
holds for \emph{all} $t \geq 0$ and for all $\omega$ from a set
$\Omega_x$ with $\mathbb{P}_x(\Omega_x)=1$. For $x=0$ and $t=0$ this
yields $\til{\beta}^i_0=0$. Next we choose $x=e_k$ (the $k$-th
unit-vector of the canonical basis in $\mathbb{R}^d$) and $t=0$. In
view of $\til{\beta}_0=0$, \eqref{e:componentdrift} implies
$\til{\beta}^i_k=0$. Hence \eqref{e:drift} holds and we are
done.
\end{proof}

The nonnegativity of $\E(X^i)$ can also be characterized completely
in terms of the parameters of $X$.
\begin{lemma}\label{l:positive}
Let $i \in \{1,\ldots,d\}$. Then $\E(X^i)$ is $\mathbb{P}_x$-a.s.\
nonnegative for all $x \in D$ if and only if
\begin{equation}\label{e:positive}
\kappa_j(\{\xi \in D: \xi_i <-1\})=0, \quad 0 \leq j \leq m.
\end{equation}
\end{lemma}

\begin{proof} Fix $x \in D$ and let $T>0$. By \cite[I.4.61]{js.87}, $\E(X^i)$ is $\mathbb{P}_x$-a.s.\ nonnegative on $[0,T]$ if and only if $\mathbb{P}_x( \,\exists\, t \in [0,T]: \Delta X_t^i < -1)=0$. By \cite[II.1.8]{js.87} and Theorem \ref{t:2.12} this in turn is equivalent to
\begin{equation}\label{e:positivecomponent}
\begin{split}
0 &= \mathbb{E}_x\left( \sum_{t \leq T} 1_{(-\infty,-1)}(\Delta X_t^i)\right)\\
 &= \mathbb{E}_x\left( 1_{(-\infty,-1)}(\xi_i)*\mu^X_T\right) \\
                                                                    &= \mathbb{E}_x\left( 1_{(-\infty,-1)}(\xi_i)*\nu_T\right) \\
                                                                    &= T\kappa_0(\{\xi \in D:\xi_i<-1\})+\sum_{j=1}^m \kappa_j(\{\xi \in D:\xi_i<-1\}) \int_0^T \mathbb{E}_x(X^j_{t-})dt.
\end{split}
\end{equation}
$\Leftarrow$: Evidently, \eqref{e:positive} implies \eqref{e:positivecomponent} for every $T$.\\
$\Rightarrow$: Since $X^j$ is nonnegative for $j=1,\ldots,m$,
\eqref{e:positivecomponent} implies that $\kappa_0(\{\xi \in D:
\xi_i<-1\})=0$ and $\kappa_j(\{\xi \in D:\xi_i<-1\})\int_0^T
\mathbb{E}_x(X^j_{t-})dt=0$ for all $x \in D$. As
in the proof of Lemma \ref{l:mloc}, it follows that $\int_0^T
\mathbb{E}_x(X^j_{t-})dt$ is strictly positive for any $x \in D$
with $x^j>0$. Hence $\kappa_j(\{\xi \in D:\xi_i<-1\})=0$, which
completes the proof.
 \end{proof}

Every positive local martingale of the form $M=\E(X^i)$ is a true
martingale for processes $X^i$ with independent increments by
\cite[Proposition 3.12]{kallsen.muhlekarbe.08b}. In general, this
does not hold true for affine processes as exemplified by
\cite[Example 3.11]{kallsen.muhlekarbe.08b}, where the following
\emph{necessary} condition is violated.
\begin{lemma}\label{l:nec}
Let $i \in \{1,\ldots,d\}$ such that $M=\E(X^i)$ is
$\mathbb{P}_x$-a.s.\ nonnegative for all $x \in D$. If $M$ is a
local $\mathbb{P}_x$-martingale for all $x \in D$, the parameters
$(\alpha^\star,\beta^\star,0,\kappa^\star)$ given by
\begin{alignat}{2}
\alpha^{\star}_j&:=\alpha_j ,&\quad  0 &\leq j \leq m,\label{e:cstar}\\
\beta^{\star}_j&:=\beta_j+\alpha_j^{\cdot i}+\int_{D \backslash \{0\}} (\xi_i h(\xi)) \kappa_j(d\xi),&\quad   0 &\leq j \leq d,\label{e:bstar}\\
\kappa^{\star}_j(d\xi) &:= (1+\xi_i)\kappa_j(d\xi),&\quad  0 &\leq j
\leq d,\label{e:nustar}
\end{alignat}
are admissible. If $M$ is a true $\mathbb{P}_x$-martingale for all
$x \in D$, the corresponding affine process
$(X,\mathbb{P}^{\star}_x)_{x \in D}$ is conservative.
\end{lemma}

\begin{proof} The first part of the assertion follows from Lemmas \ref{l:mloc} and \ref{l:positive} as in the proof of \cite[Lemma 3.5]{kallsen.muhlekarbe.08b}. Let $M$ be a true martingale for all $x \in D$. Then for every $x \in D$, e.g.\ \cite{cherny.02} shows that there exists a probability measure $\mathbb{P}^M_x \stackrel{\mathrm{loc}}{\ll} \mathbb{P}_x$ on $(\mathbb{D}^d,\scr{D}^d,(\scr{D}^d_t))$ with density process $M$. Then the Girsanov-Jacod-Memin theorem as in \cite[Lemma 5.1]{kallsen.03}  yields that $X$ admits affine $\mathbb{P}^M_x$-characteristics as in \eqref{e:b}-\eqref{e:nu} with $(\alpha,\beta,0,\kappa)$ replaced by $(\alpha^{\star},\beta^{\star},0,\kappa^{\star})$. Since $\mathbb{P}^M_x |_{\scr{D}_0} = \mathbb{P}_x |_{\scr{D}_0}$ implies $\mathbb{P}_x^M(X_0=x)=1$, we have $\mathbb{P}^M_x=\mathbb{P}^{\star}_x$ by Theorem \ref{t:2.12} . In particular, the transition function $p_t^{\star}(x,d\xi)$ of $(X,\mathbb{P}^{\star}_x)_{x \in D}$ satisfies $1=\mathbb
 {P}_x^M(X_t \in D)=\mathbb{P}_x^{\star}(X_t \in D)=p^{\star}_t(x,D)$, which completes the proof.\end{proof}
If $M=\E(X^i)$ is only a local martingale, the affine process
$(X,\mathbb{P}_x^{\star})_{x \in D}$ does not necessarily have to be
conservative (see \cite[Example 3.11]{kallsen.muhlekarbe.08b}). A
careful inspection of the proof of \cite[Theorem
3.1]{kallsen.muhlekarbe.08b} reveals that conservativeness of
$(X,\mathbb{P}^{\star}_x)_{x \in D}$ is also a \emph{sufficient}
condition for $M$ to be a martingale. Combined with Lemma
\ref{l:mloc} and Theorem \ref{th: char can cons} this in turn allows
us to provide the following deterministic necessary and sufficient
conditions for the martingale property of $M$ in terms of the
parameters of $X$.
\begin{theorem}\label{t:main}
Let $i \in \{1,\ldots,d\}$ such that $\E(X^i)$ is
$\mathbb{P}_x$-a.s.\ nonnegative for all $x \in D$. Then we have
equivalence between:
\begin{enumerate}
\item $\E(X^i)$ is a true $\mathbb{P}_x$-martingale for all $x \in D$.\label{equiv:1}
\item $\E(X^i)$ is a local $\mathbb{P}_x$-martingale for all $x \in D$ and the affine process corresponding to the admissible parameters $(\alpha^{\star},\beta^{\star},0,\kappa^{\star})$ given by  \eqref{e:cstar}-\eqref{e:nustar} is conservative.\label{equiv:2}
\item \eqref{e:integrable} and \eqref{e:drift} hold and $g=0$ is the only $\rr^m_{-}$-valued local solution of
\begin{equation}\label{e:eberhard}
\partial_t g(t)=R^{\star}_{\mathcal I}(g(t),0), \quad g(0)=0,
\end{equation}
 where $R^\star$ is given by \eqref{e:R} with $(\alpha^{\star},\beta^{\star},0,\kappa^{\star})$ instead of $(\alpha,\beta,\gamma,\kappa)$.\label{equiv:3}
\end{enumerate}
\end{theorem}

\begin{proof}  \ref{equiv:1} $\Rightarrow$ \ref{equiv:2}: This is shown in Lemma \ref{l:nec}.\\
\ref{equiv:2} $\Rightarrow$ \ref{equiv:3}: This follows from Lemma
\ref{l:mloc} and Theorem \ref{th: char can cons}.\\\ref{equiv:3}
$\Rightarrow$ \ref{equiv:1}: By \eqref{e:integrable},
\eqref{e:drift} and Lemma \ref{l:positive}, Assumptions 1-3 of
\cite[Theorem 3.1]{kallsen.muhlekarbe.08b} are satisfied. Since we
consider time-homogeneous parameters here, Condition 4 of
\cite[Theorem 3.1]{kallsen.muhlekarbe.08b} also follows immediately
from \eqref{e:integrable}. The final Condition 5 of \cite[Theorem
3.1]{kallsen.muhlekarbe.08b} is only needed in \cite[Lemma
3.5]{kallsen.muhlekarbe.08b} to ensure that a semimartingale with
affine characteristics relative to
$(\alpha^{\star},\beta^\star,0,\kappa^\star)$ exists. In view of the
first part of Lemma \ref{l:nec}, Theorem \ref{th: char can cons} and
Theorem \ref{t:2.12} it can therefore be replaced by requiring that
$0$ is the unique $\rr^m_{-}$-valued solution to \eqref{e:eberhard}.
The proof of \cite[Theorem 3.1]{kallsen.muhlekarbe.08b} can then be
carried through unchanged. \end{proof}

\begin{remark}\label{remarkmart}\rm
\begin{enumerate}
\item \label{Remark1} In view of \cite[Lemma 2.7]{kallsen.muhlekarbe.08b},  $\til{M}:=\exp(X^i)$ can be written as  $\til{M}=\exp(X^i_0)\E(\til{X}^i)$ for the $d+1$-th component of the $\rp^m \times \rr^{n+1}$-valued affine process $(X,\til{X}^i)$ corresponding to the admissible parameters $(\til{\alpha},\til{\beta},0,\til{\kappa})$ given by $(\til{\alpha}_{d+1},\til{\beta}_{d+1},\til{\kappa}_{d+1})=(0,0,0)$ and
$$\qquad \quad  (\til{\alpha}_j,\til{\beta}_j,\til{\kappa}_j(G)):=\left(\begin{pmatrix} \alpha_j & \alpha_j^{\cdot i} \\ \alpha_j^{i \cdot} & \alpha_j^{ii} \end{pmatrix} , \begin{pmatrix} \beta_j \\ \til{\beta}^{d+1}_j \end{pmatrix},\int_{D \backslash \{0\}} 1_G(\xi,e^{\xi_i}-1)\kappa_j(d\xi)\right)$$
for $G \in \B^{d+1}$, $j=0,\ldots,d$, and
\begin{equation*}
\til{\beta}^{d+1}_j=\beta_j^i+\frac{1}{2}\alpha_j^{ii}+\int_{D
\backslash \{0\}} (h_i(e^{\xi_i}-1)-h_i(\xi))\kappa_j(d\xi).
\end{equation*}
This allows to apply Theorem \ref{t:main} in this situation as well.
\item Theorem \ref{t:main} is stated for the stochastic exponential $\E(X^i)$ of $X^i$, that is, the projection of $X$ to the $i$-th component. It can, however, also be applied to the stochastic exponential $\E(A(X))$ of a general affine functional $A:D \to \mathbb{R}: x \mapsto p+Px$, where $p \in \mathbb{R}$ and $P \in \mathbb{R}^d$. To see this, note that it follows from It\^o's formula and Theorem \ref{t:2.12} that the $\mathbb{R}^m_+ \times \mathbb{R}^{n+1}$-valued process $Y=(X,A(X))$ is affine with admissible parameters $(\til{\alpha},\til{\beta},0,\til{\kappa})$ given by $(\til{\alpha}_{d+1},\til{\beta}_{d+1},\til{\kappa}_{d+1})=(0,0,0)$ and
$$
\qquad \qquad \til{\alpha}_j=\begin{pmatrix} \alpha_j & \alpha_j P \\  P^{\top} \alpha_j & P^{\top} \alpha_j P \end{pmatrix}, \quad \til{\beta}_j = \begin{pmatrix} \beta_j \\  P^{\top} \beta_j + \int (h(P^{\top}x)-P^{\top}h(x))\kappa_j(dx) \end{pmatrix}, 
$$
as well as 
$$\til{\kappa}_j(G)=\int_{D\backslash \{0\}} 1_G(x,P^{\top}x)\kappa_j(dx) \quad \forall G \in \scr{B}^{d+1},$$
for $j=0,\ldots,d$. Therefore one can simply apply Theorem \ref{t:main} to $\E(Y^{d+1})$.
\item \label{expremark2} Conservativeness of $(X,\mathbb{P}^\star_x)_{x \in D}$ and uniqueness for solutions of ODE \eqref{e:eberhard1} can be ensured by requiring the moment condition \eqref{e:sufficient}
    for $\kappa_j^\star$. The implication \ref{equiv:3} $\Rightarrow$ \ref{equiv:1} in Theorem \ref{t:main} therefore leads to the easy-to-check sufficient criterion  \cite[Corollary 3.9]{kallsen.muhlekarbe.08b} for the martingale property of $M$.
\item \label{expremark3} By Remark \ref{consremark} \ref{expremark3cons} we know that in the case $m=1$, $(X,\mathbb{P}^\star_x)_{x \in D}$ is conservative if and only if either \eqref{e:sufficient} holds for $\kappa_j^\star$ or equation \eqref{eq osgood} holds for $R_{1}^{\star}$. Together with Remark \ref{Remark1}, this leads to the necessary and sufficient condition for the martingale property of ordinary exponentials $\exp(X^i)$ obtained in \cite[Theorem 2.5]{kellerressel.09}.
\end{enumerate}
\end{remark}

We conclude by providing an example of an exponentially affine local
martingale for which the sufficient conditions from
\cite{kallsen.muhlekarbe.08b} cannot be applied.  Our main Theorem
\ref{t:main}, however, shows that is indeed a true martingale. This process is based on the one in Example \ref{crucial cons example} and therefore again somewhat artificial. Various extensions are possible, but we again restrict ourselves to the simplest possible specification here.

\begin{example}\rm
Consider the $\mathbb{R}_+ \times \mathbb{R}$-valued affine process
$(X^1,X^2)$ corresponding to the admissible parameters
$$\alpha=(0,0,0), \quad \beta=(0,\beta_1,0), \quad \gamma=(0,0,0),\quad  \kappa=(0,\kappa_1,0),$$
where
$$\begin{pmatrix} \beta^1_1 \\ \beta^2_1 \end{pmatrix}=\begin{pmatrix} \sum_{n=1}^\infty  \frac{1}{(1+n)n^2} \\ \sum_{n=1}^\infty \frac{1-n}{(1+n)n^2}\end{pmatrix} \quad \mbox{and} \quad \kappa_1=\sum_{n=1}^\infty \frac{\delta_{(n,n)}}{(1+n)n^{2}},$$
for the Dirac measures $\delta_{(n,n)}$ supported by $\{(n,n)\}$, $n \in \mathbb{N}$. Since
$X^2$ has only positive jumps, $\E(X^2)$ is positive. Moreover, it
is a local martingale by Lemma \ref{l:mloc}, because
$$\int_{\{|\xi_2|>1\}}|\xi_2|\kappa_1(d\xi)= \sum_{n=1}^\infty \frac{1}{(1+n)n}<\infty$$
and $\beta^2_1+\int_0^\infty (\xi_2-h_2(\xi_2))\kappa_1(d\xi)=0$. Note
that \cite[Corollary 3.9]{kallsen.muhlekarbe.08b} is not applicable, because
$$\int_{\{|\xi_2|>1\}}\xi_1(1+\xi_2) \kappa_1(d\xi)= \sum_{n=1}^\infty \frac{1}{n}=\infty.$$
However, by Theorem \ref{t:main} and Remark \ref{remarkmart}(iii),
$\E(X^2)$ is a true martingale, since we have shown in Example
\ref{crucial cons example} that \eqref{eq osgood} is satisfied for
$$R^\star_1(u_1,0)=\sum_{n=1}^\infty \frac{e^{u_1 n}-1}{n^2}.$$
\end{example}

\begin{appendix}
\section{ODE comparison results in non-Lipschitz setting}\label{sec: ODEs}
Let $C$ be a closed convex proper cone with nonempty interior
$C^\circ$ in a normed vector space $(E,\|\,\,\|)$. The partial order
induced by $C$ is denoted by $\preceqq$. For $x,y\in E$, we write
$x\ll y$ if $y-x\in C^\circ$. We denote by $C^*$ the dual cone of
$C$. Let $D_g$ be a set in $E$. A function $g\colon D_g\rightarrow
E$ is called \emph{quasimonotone increasing}, in short \emph{qmi},
if for all $l\in C^*$, and $x,y\in D_g$
\[
(x\preceqq y,\,\,l(x)=l(y))\Rightarrow(l(g(x))\leq l(g(y))).
\]
The next lemma is a special case of Volkmann's result \cite[Satz
1]{Volkmann}.
\begin{lemma}\label{th: Volkmann}
Let $0<T\leq \infty$, $D_f\subset E$, and $f\colon [0,T)\times
D_f\rightarrow E$ be such that $f(t,\cdot)$ is qmi on $D_f$ for all
$t\in [0,T)$. Let $\zeta,\eta:[0,T)\rightarrow D_f$ be curves that
are continuous on $[0,T)$ and differentiable on $(0,T)$. Suppose
$\zeta(0)\gg \eta(0)$ and $\dot {\zeta}(t)-f(t,\zeta(t))\gg
\dot{\eta}(t)-f(t,\eta(t))$ for all $t\in (0,T)$. Then $\zeta(t)\gg
\eta(t)$ for all $t\in[0,T)$.
\end{lemma}
A function $g:[0,T)\times D_g\rightarrow E$ is called \emph{locally
Lipschitz}, if for all $0<t<T$ and for all compact sets $K\subset
D_g$ we have
\[
L_{t,K}(g):=\sup_{0<\tau<t,\ x,y \in K: x \neq
y}\frac{\|g(\tau,x)-g(\tau,y)\|}{\|x-y\|}<\infty
\]
where $L_{t,K}(g)$ is usually called the Lipschitz constant.

We now use Lemma \ref{th: Volkmann} to prove the following general
comparison result.
\begin{proposition}\label{prop: essential comparison}
Let $T$, $D_f$, and $f$ be as in Lemma~$\ref{th: Volkmann}$.
Suppose, moreover, that $D_f$ has a nonempty interior and $f$ is
locally Lipschitz on $[0,T)\times D_f^\circ$. Let
$\zeta,\eta:[0,T)\rightarrow D_f$ be curves that are continuous on
$[0,T)$, differentiable on $(0,T)$, and satisfy the conditions
\begin{enumerate}
\item $\eta(t)\in D_f^\circ$

\item $\dot {\zeta}(t)-f(t,\zeta(t))\succeqq \dot{\eta}(t)-f(t,\eta(t))$

\item $\zeta(0)\succeqq \eta(0)$
\end{enumerate}
for all $t\in [0,T)$. Then $\zeta(t)\succeqq\eta(t)$ for all $t\in
[0,T)$.
\end{proposition}
\begin{proof}
Fix $t_0\in [0,T)$. Since $\eta$ is continuous, the image $S$ of the
segment $[0,t_0]$ under the map $\eta$ is a compact subset of
$D_f^\circ$. Let $\delta>0$ be such that the closed
$\delta$-neighborhood $S_\delta$ of $S$ is contained in $D_f^\circ$.
By the local Lipschitz continuity of $f$ on $D_f^\circ$, there
exists a constant $L>0$ such that
\begin{equation}\label{Lipschitz}
\|f(t,x)-f(t,y)\|\leq L\|x-y\|
\end{equation}
for any $t\in [0,t_0]$ and $x,y\in S_\delta$. Let $c\in C^\circ$ be
such that $\|c\|=1$ and let $d_c$ denote the distance from $c$ to
the boundary $\partial C$ of $C$. For $\varepsilon>0$, we set
$h_\varepsilon(t):=\varepsilon e^{2Lt/d_c} c$. If $\varepsilon \leq
e^{-2Lt_0/d_c}\delta$, then $\eta(t)-h_\varepsilon(t)\in S_{\delta}$
for any $t\in [0,t_0]$, and (\ref{Lipschitz}) gives
\begin{equation}\label{norm_estimate}
\|f(t,\eta(t)-h_\varepsilon(t))-f(t,\eta(t))\|\leq
L\|h_{\varepsilon}(t)\|,\quad t\in [0,t_0].
\end{equation}
Since $C$ is a cone, the distance from $Lh_\varepsilon(t)/d_c$ to
$\partial C$ is equal to $L\varepsilon
e^{2Lt/d_c}=L\|h_\varepsilon(t)\|$. In view
of~(\ref{norm_estimate}), it follows that
\[
Lh_\varepsilon(t)/d_c\succeqq
f(t,\eta(t)-h_\varepsilon(t))-f(t,\eta(t))
\]
and hence
\begin{equation}\label{eq: comp1}
-\dot h_\varepsilon(t)=-2L h_\varepsilon(t)/d_c\ll
f(t,\eta(t)-h_\varepsilon(t))-f(t,\eta(t)),\quad t\in [0,t_0],
\end{equation}
for $\varepsilon$ small enough. This implies that
\begin{equation*}
\dot{\zeta}(t)-f(t,\zeta(t))\succeqq \dot{\eta}(t)-f(t,\eta(t))\gg
\dot{\eta}(t)-\dot h_\varepsilon(t)-f(t,\eta(t)+h_\varepsilon(t)).
\end{equation*}
Applying Lemma \ref{th: Volkmann} to the functions $\zeta(t)$ and
$\eta(t)+h_\varepsilon(t)$ yields $\zeta(t)\gg
\eta(t)+h_\varepsilon(t)$, for all $t\in [0, t_0]$. Now letting
$\varepsilon\rightarrow 0$ yields the required inequality for all
$t\in [0,t_0]$. This proves the assertion, because $t_0< T$ can be
chosen arbitrarily.
\end{proof}

If we consider the differential equation
\begin{equation}\label{eq: ODE general}
\dot {\xi}=f(t,\xi(t)),\quad \xi(0)=u\in D_f,
\end{equation}
Proposition \ref{prop: essential comparison} allows the following
immediate conclusion, which is the key tool for proving Proposition
\ref{prop: extremality} and in turn Theorem \ref{th: char can cons}.
\begin{corollary}\label{th1}
Let $T$, $D_f$ and $f$ be as in Lemma \ref{prop: essential
comparison}. Suppose further that equation \eqref{eq: ODE general}
gives rise to a global solution $\psi^\circ(t,u)\colon\,\mathbb
R_+\times D^\circ_f\rightarrow D_f^\circ$. Let $u_2\in D_f^\circ$
and let $\xi\colon [0,T)\rightarrow D_f$ be a solution of \eqref{eq:
ODE general} such that $\xi(0)=u_1\succeqq u_2$. Then
$\xi(t)\succeqq \psi^\circ(t,u_2)$, for all $t\in [0,T)$.
\end{corollary}

\end{appendix}

\providecommand{\bysame}{\leavevmode\hbox
to3em{\hrulefill}\thinspace}
\providecommand{\MR}{\relax\ifhmode\unskip\space\fi MR }
\providecommand{\MRhref}[2]{%
  \href{http://www.ams.org/mathscinet-getitem?mr=#1}{#2}
} \providecommand{\href}[2]{#2}

\end{document}